\documentclass[number,citesort,dvips]{arxbj}
\usepackage[vtexbibl]{}

\aid{0}
\volume{16}
\issue{4}
\pubyear{2010}
\firstpage{1137}
\lastpage{1163}
\doi{10.3150/09-BEJ239}

\makeatletter

\newtheorem{theorem}{Theorem}
\newremark{remark}{Remark}
\newtheorem{proposition}{Proposition}

\newtheorem{corollary}{Corollary}

\newtheorem{lemma}{Lemma}

\makeatother

\begin{document}
\begin{frontmatter}

\title{Adaptive estimation of a distribution function and its density
in sup-norm loss by wavelet and spline projections}
\runtitle{Adaptive estimation}

\begin{aug}
\author[a]{\fnms{Evarist} \snm{Gin\'{e}}\thanksref{a}\ead[label=e1]{gine@math.uconn.edu}\corref{}}
\and
\author[b]{\fnms{Richard} \snm{Nickl}\thanksref{b}\ead[label=e2]{nickl@statslab.cam.ac.uk}}
\runauthor{E. Gin\'{e} and R. Nickl}
\address[a]{Department of Mathematics, University of Connecticut,
Storrs, CT 06269-3009, USA.\\ \printead{e1}}
\address[b]{Statistical Laboratory, Department of Pure Mathematics and
Mathematical Statistics,
University of Cambridge, Wilberforce Road, Cambridge CB3 0WB, UK.
\printead{e2}}
\end{aug}

% HISTORY:
\received{\smonth{5} \syear{2008}}
\revised{\smonth{8} \syear{2009}}

% ABSTRACT
%
\begin{abstract}
Given an i.i.d.~sample from a distribution $F$ on $\mathbb R$ with
uniformly continuous density $p_0$, purely
data-driven estimators are constructed that efficiently estimate $F$ in
sup-norm loss and simultaneously estimate
$p_0$ at the best possible rate of convergence over H\"{o}lder balls,
also in sup-norm loss. The estimators are
obtained by applying a model selection procedure close to Lepski's
method with random thresholds to projections
of the empirical measure onto spaces spanned by wavelets or
$B$-splines. The random thresholds are based on
suprema of Rademacher processes indexed by wavelet or spline projection
kernels. This requires Bernstein-type
analogs of the inequalities in Koltchinskii [\textit
{Ann.~Statist.}~\textbf{34} (2006) 2593--2656] for the deviation
of suprema of empirical processes from their Rademacher symmetrizations.
\end{abstract}

% KEYWORDS
%

\end{frontmatter}

%s1 ###
\section{Introduction}\label{sec1}

If $X_1,\ldots,X_n$ are i.i.d.~with unknown distribution function $F$ on
$\mathbb R$, then classical results of mathematical statistics
establish optimality of the empirical distribution function $F_n$ as an
estimator of $F$. That is to say, if we assume no a priori knowledge
whatsoever on $F$ and equip the set of all probability distribution
functions with some natural loss function such as sup-norm loss,
 then~$F_n$ is asymptotically sharp minimax for estimating $F$. (The same is
true even if more is known about~$F$, for instance, if $F$ is known to
have a uniformly continuous density.) However, this does not preclude
the existence of other estimators that are also asymptotically minimax
for estimating $F$ in sup-norm loss, but which improve upon $F_n$ in
other respects. What we have in mind is a purely data-driven estimator
that is efficient for $F$, but, at the same time, also estimates the
density $f$ of $F$ at the best rate of convergence in some relevant
loss function over some prescribed classes of densities. More
precisely, our goal in the present article is to construct estimators
that satisfy the functional central limit theorem (CLT) for the
distribution function \textit{and} which adapt to the unknown
smoothness of the density in \textit{sup-norm loss}. Whereas this
article is concerned with the mathematical problem of the existence and
construction of such estimators, it does not deal with the practical
implementation of estimation procedures.

To achieve adaptation, one can opt for several approaches, all of which
are related. Among them, we mention the penalization method of Barron,
Birg\'{e} and Massart \cite{1999Barron}, wavelet threshholding \cite{1996Donoho} and Lepski's \cite{1991Lepski} method. Our choice
for the goal at hand consists of using Lepski's method, with random
thresholds, applied to wavelet and spline projection estimators of a density.

The linear estimators underlying our procedure are projections of the
empirical measure onto spaces spanned by wavelets, and wavelet theory
is central to some of the derivations of this article. The wavelets
most commonly used in statistics are those that are compactly supported
(for example, Daubechies wavelets), and our results readily apply to
these. However, for computational and other purposes, projections onto
spline spaces are also interesting candidates for the estimators.
Density estimators obtained by projecting the empirical measure onto
Schoenberg spaces spanned by $B$-splines were studied by Huang and
Studden \cite{1993Huang}. As is well known in wavelet theory, the Schoenberg
spline spaces with equally spaced knots have an orthonormal basis
consisting of the Battle--Lemari\'{e} wavelets so that the spline
projection estimator is, in fact, exactly equal to the wavelet
estimator based on Battle--Lemari\'{e} wavelets. These wavelets do not
have compact support, but they are exponentially localized. Although we
cannot, in general, handle exponentially decaying wavelets, we can
still work with Battle--Lemari\'e wavelets because the $B$-spline
expansion of the projections allows us to show that the relevant
classes of functions are of Vapnik--Chervonenkis type so that empirical
process techniques can be applied. In particular, the adaptive
estimators we devise in Theorem \ref{cdf} may be based either on spline
projections or on compactly supported wavelets. In the process of
proving the main theorem, we also provide new asymptotic results for
spline projection density estimators similar to those for wavelet
estimators in \cite{2009bGin}.

We need to use Talagrand's exponential inequality with sharp constants
\cite{2003Bousquet,2005Klein} in the proofs, but to do this,
we have to estimate the expectation of suprema of certain empirical
processes that appear in the centering of Talagrand's inequality. The
use of entropy-based moment inequalities for empirical processes
typically results in too conservative constants (for example, in \cite{2009aGin}). In order to remedy this problem, we adapt
recent ideas due to Koltchinskii \cite{2001Koltchinskii,2006Koltchinskii} and Bartlett, Boucheron
and Lugosi \cite{2002Bartlett} to density estimation: the entropy-based moment
bounds are replaced by the sup-norm of the associated Rademacher
averages, which are, with high probability, better estimates of the
expected value of the supremum of the empirical process. We derive a
Bernstein-type analog of an exponential inequality in \cite{2006Koltchinskii} that shows how the supremum of an empirical process deviates
from the supremum of the associated Rademacher processes. This
Bernstein-type version allows one to use partial knowledge of the
variance of the empirical processes involved, which is crucial for
applications in our context of adaptive density estimation. Moreover,
we show that one can use, instead of the supremum of the Rademacher
process, its conditional expectation given the data.

Adaptive estimation in sup-norm loss is a relatively recent subject. We
should mention the results in Tsybakov \cite{1998Tsybakov}, Golubev, Lepski and
Levit \cite{2001Golubev} -- who only considered Sobolev-type smoothness conditions
-- and  \cite{2009Goldenshluger}. All of these results were
obtained in the Gaussian white noise model. If one is interested in
adapting to a H\"{o}lder-continuous density in sup-norm loss in the
i.i.d.~density model on $\mathbb R$, this simplifying Gaussian
structure is not available and novel techniques are needed. In the
i.i.d.~density model on $\mathbb R$, a direct `competitor' to the
estimators constructed in this article is the hard thresholding wavelet
density estimator introduced in \cite{1996Donoho}: as proved in \cite{2009bGin}, its distribution function satisfies the functional
CLT and it is adaptive in the sup-norm over H\"older balls; however,
the proofs there seem to require the additional assumption that $\mathrm{d}F$ integrates
$|x|^\delta$ for some $\delta>0$, and the constants appearing in the
threshold and the risk become quite large for $\delta$ small. The
results in the present article hold under no moment condition whatsoever.

%s2 ###
\section{Wavelet expansions and estimators}\label{sec2}

We start with some basic notation. If $(S, \mathcal S)$ is a measurable
space, then for Borel-measurable functions $h\dvtx S \to\mathbb{R}$ and
Borel measures $\mu$ on $S$, we set $\mu h:=\int_S h \,\mathrm{d}\mu$. We will
denote by $L^p(Q):= L^p(S,Q)$, $1 \le p \le\infty$, the usual Lebesgue
spaces on $S$ with respect to~a Borel measure $Q$, and if $Q$ is
Lebesgue measure on $S = \mathbb R$, then we simply denote this space
by $L^p(\mathbb R)$, and its norm by $\|\cdot\|_p$, if $p < \infty$. We
will use $\|h\|_\infty$ to denote $\sup_{x \in\mathbb R} |h(x)|$ for
$h\dvtx \mathbb R \to\mathbb R$. For $s \in\mathbb N$, denote by $\mathsf
{C}^{s}(\mathbb R)$ the spaces of functions $f\dvtx \mathbb R \to\mathbb R$
that are $s$-times differentiable with bounded uniformly continuous
$D^rf$, $0<r\leq s$,
equipped with the norm $ \| f\|_{s,\infty}=\sum_{0\leq\alpha\leq
s}\Vert D^{\alpha}f\Vert_{\infty},$ with the convention that $D^0 =:
id$ and that $\mathsf C (\mathbb R) := \mathsf C^0(\mathbb R)$ is then
the space of bounded uniformly continuous functions. For non-integer
${s>0}$ and $[ s] $ the integer part of $s$, set
\begin{eqnarray*}
\mathsf{C}^{s}(\mathbb R)=\biggl\{  f\in\mathsf{C}^{[s]}(\mathbb{R}
)\dvtx \Vert f\Vert_{s,\infty}:=\sum_{0\leq\alpha\leq[ s]}\Vert
D^{\alpha}f\Vert_{\infty}+ \sup_{x \ne y} \frac{\vert
D^{[s]}f(x)-D^{[s]}f(y)\vert}{\vert
x-y\vert^{s-[ s] }}<\infty  \biggr\} .
\end{eqnarray*}

%s2.1 ###
\subsection{Multiresolution analysis and wavelet bases}

We recall here a few well-known facts about wavelet expansions; see,
for example, Sections 8 and 9 in \cite{1998Har}. Let $\phi\in L^2(\mathbb R)$ be a scaling
function, that is, $\phi$ is such that $\{\phi(\cdot-k)\dvtx k \in\mathbb
Z\}$ is an orthonormal system in $L^2 (\mathbb R)$ and, moreover, the
linear spaces $V_0 = \{f(x)= \sum_k c_k \phi(x-k) \dvtx \{c_k\}_{k \in
\mathbb Z} \in\ell^2 \}$, $V_1 = \{h(x)=f(2x) \dvtx f \in V_0\},\ldots,V_j
=\{h(x)=f(2^jx) \dvtx f \in V_0\},\ldots$ are nested ($V_{j-1} \subseteq
V_{j}$ for $j \in\mathbb N$) and their union is dense in $L^2(\mathbb R)$.
In the case where $\phi$ is a bounded function that decays
exponentially at infinity (that is, $|\phi(x)| \le C \mathrm{e}^{-\gamma|x|}$
for some $C, \gamma>0$) -- which we assume for the rest of this
subsection -- the kernel of the projection onto the space $V_j$ has
certain properties. First, the series
%
%e1 ###
\begin{equation}\label{kernel0}
K(y,x) := K(\phi,y,x) = \sum_{k \in\mathbb Z} \phi(y-k) \phi(x-k)
\end{equation}
converges pointwise and we set $K_j(y,x) := 2^j K(2^jy,2^jx), j \in
\mathbb N \cup\{0\}$. Furthermore, we have
%
%e2 ###
\begin{equation} \label{major}
|K(y,x)| \leq\Phi(|y-x|)    \quad\mbox{and} \quad    \sup_{x \in\mathbb R}
\sum_k \bigl|\phi(x-k)\bigr| < \infty,
\end{equation}
where $\Phi\dvtx \mathbb R \to\mathbb R^+$ is bounded and has exponential
decay (cf. Lemma 8.6 in \cite{1996Talagrand}). For any $j$ fixed, if $f\in
L^p(\mathbb R)$, $1\le p\le\infty$, then the series
\[
K_j(f)(y) := \int K_j(x,y) f(x)\,\mathrm{d}x = \sum_{k \in\mathbb Z} 2^j\phi
(2^jy-k)\int\phi(2^jx-k)f(x)\,\mathrm{d}x,\qquad    y \in\mathbb R ,
\]
converges pointwise and, for $f \in L^2(\mathbb R)$, $K_j(f)$ coincides
with the orthogonal projection $\pi_j \dvtx L^2(\mathbb R) \to V_j$ of $f$
onto $V_j$. For $f \in L^1(\mathbb R)$, which is the main case in this
article, the convergence of the series in fact takes place in
$L^p(\mathbb R)$, $1\le p\le\infty$. This still holds true if $f(x)\,\mathrm{d}x$
is replaced by $\mathrm{d}\mu(x)$, where $\mu$ is any finite signed measure.
If, now, $\phi$ is a scaling function and $\psi$ the associated mother
wavelet so that $\{\phi(\cdot-k), 2^{l/2} \psi(2^l(\cdot)-k)\dvtx k \in
\mathbb Z, l \in\mathbb N \}$ is an orthonormal basis of $L^2(\mathbb
R)$, then any $f \in L^p(\mathbb R)$ admits the formal expansion
%
%e3 ###
\begin{equation} \label{wavser}
f(y) = \sum_k \alpha_k(f) \phi(y-k) + \sum_{l=0}^\infty\sum_k \beta
_{lk} (f) \psi_{lk}(y) ,
\end{equation}
where $\psi_{lk}(y) = 2^{l/2} \psi(2^ly-k)$, $\alpha_k(f) = \int f(x)
\phi(x-k)\,\mathrm{d}x$, $\beta_{lk}(f) = \int f(x) \psi_{lk}(x)\,\mathrm{d}x$. Since
$(K_{l+1}-K_l)f = \sum_k \beta_{lk}(f) \psi_{lk}$, the partial sums of
the series (\ref{wavser}) are in fact given by
%
%e4 ###
\begin{equation}\label{wavexp}
K_j(f)(y) = \sum_k \alpha_k(f) \phi(y-k) + \sum_{l=0}^{j-1} \sum_k \beta
_{lk}(f) \psi_{lk}(y)
\end{equation}
and if $\phi, \psi$ are bounded and have exponential decay, then
convergence of the series (\ref{wavexp}) holds pointwise; it also holds
in $ L^p(\mathbb R)$, $1 \leq p \leq\infty$, if $f \in L^1 (\mathbb
R)$ or if $f$ is replaced by a finite signed measure. Now, using these
facts, one can furthermore show that the wavelet series (\ref{wavser})
converges in $ L^p(\mathbb R)$, $p<\infty$, for $f \in L^p(\mathbb R)$
and we also note that if $p_0$ is a uniformly continuous density, then
its wavelet series converges uniformly.

%s2.2 ###
\subsection{Density estimation using wavelet and spline projection kernels}

Let $X_1,\ldots,X_n$ be i.i.d. random variables with common law $P$ and
density $p_0$ on $\mathbb R$, and denote by $P_n = \frac{1}{n}\sum
_{i=1}^n \delta_{X_i}$ the associated empirical measure. A natural
first step is to estimate the projection $K_j(p_0)$ of $p_0$ onto $V_j$ by
%
%e5 ###
\begin{equation}\label{est0}
p_n(y) := p_n (y,j) =\frac{1}{n} \sum_{i=1}^n K_j (y, X_i) = \sum_k
\hat\alpha_k \phi(y-k) + \sum_{l=0}^{j-1} \sum_k \hat\beta_{lk} \psi
_{lk}(y) , \qquad  y \in\mathbb R,\hspace*{4pt}
\end{equation}
where $K$ is as in (\ref{kernel0}), $j \in\mathbb N$, and where $\hat
\alpha_k = \int\phi(x-k)\,\mathrm{d}P_n(x)$, $\hat\beta_{lk} = \int\psi
_{lk}(x)\,\mathrm{d}P_n(x)$ are the empirical wavelet coefficients. We note that
for $\phi$, $\psi$ compactly supported (for example, Daubechies
wavelets), there are only finitely many $k$'s for which these
coefficients are non-zero. This estimator was first studied by
Kerkyacharian and Picard \cite{1992Kerkyacharian} for compactly supported wavelets.

If the wavelets $\phi$ and $\psi$ do not have compact support, it may
be impossible to compute the estimator exactly since the sums over $k$
consist of infinitely many summands. However, in the special case of
the Battle--Lemari\'{e} family $\phi_r, r \geq1$ (see, for example,
Section 6.1 in \cite{1998Har}) -- which is a class of non-compactly
supported but exponentially decaying wavelets -- the estimator has a
simple form in terms of splines: the associated spaces $V_{j,r} = \{\sum
_k c_k 2^{j/2}\phi_r(2^j(\cdot)-k) \dvtx \sum_k c^2_k < \infty\}$ are, in
fact, equal to the \textit{Schoenberg spaces} generated by the Riesz
basis of $B$-splines of order $r$ so that the sum in (\ref{est0}) can
be computed by
%
%e6 ###
\begin{equation} \label{est}
p_n(y,j) := \frac{1}{n} \sum_{i=1}^n \kappa_j (y, X_i) = \frac{2^j}{n}
\sum_{i=1}^n \sum_{k} \sum_{l} b_{kl} N_{j,k,r}(X_i) N_{j,l,r}(y)
, \qquad    y \in\mathbb R,
\end{equation}
where the $N_{j,k,r}$ are (suitably translated and dilated) $B$-splines
of order $r$, the kernel $\kappa$ is as in~(\ref{prok}) below and the
$b_{kl}$'s are the entries of the inverse of the matrix defined in (\ref
{toep}) below. An exact derivation of this spline projection, its
wavelet representation and detailed definitions are given in Section
\ref{splines}. It turns out that for every sample point $X_i$ and for
every $y$, each of the last two sums extends over only $r$ terms. We
should note that this `spline projection' estimator was first studied
(outside the wavelet setting) by Huang and Studden \cite{1993Huang}, who derived
pointwise rates of convergence; see also \cite{1999Huang}, where some
comparison between Daubechies and spline wavelets can be found.

In the course of proving the main theorem of this article, we will
derive some basic results for the linear spline projection estimator
(\ref{est}), which we now state. For classical kernel estimators,
results similar to those that follow were obtained in
\cite{2000Deheuvels,2002Gin,2009aGin}, and for wavelet estimators based on compactly supported
wavelets, this was done in \cite{2009bGin}.

\begin{theorem} \label{linear}
Suppose that $P$ has a bounded density $p_0$. Assume that $j_n \to
\infty$, $n /(j_n 2^{j_n}) \to\infty$, $j_n /\log\log n \to\infty$
and $j_{2n} - j_n \le\tau$ for some $\tau$ positive. Let $p_n(y)=
p_n(y, j_n)$ be the estimator from (\ref{est}) for some $r \geq1$. Then
\begin{eqnarray*}
\limsup_n \sqrt{\frac{n}{2^{j_n}j_n}} \sup_{y \in\mathbb R} |p_n(y) -
Ep_n(y)| = C  \qquad  \mbox{a.s.}
\end{eqnarray*}
and, for $1 \leq p < \infty$,
\begin{eqnarray*}
\sup_n \sqrt{\frac{n} {2^{j_n}j_n}} \Bigl(E\sup_{y \in\mathbb R} |p_n(y) -
Ep_n(y)|^p\Bigr)^{1/p} \le C' ,
\end{eqnarray*}
where $C$ and $C'$ depend only on $\|p_0\|_\infty$ and on $r, p, \tau$  ... and on $r$, $p$, $\tau$. Moreover, if $p_0\in C^t(\mathbb R)$, then
\[
\sup_{y\in\mathbb R}|p_n(y)-p_0(y)|=O\Biggl(\sqrt{\frac{2^{j_n}j_n}{n}}+2^{-tj_n}\Biggr)\qquad \mbox{a.s.  and in }
L^p(P).
\]
\end{theorem}

For rates of convergence in probability, the conditions on $j_n$ can be
weakened (see Proposition~\ref{ineq} below). The last bound in this theorem gives, for $p_0 \in\mathsf C^t(\mathbb R)$ with
$t \le r$ and $2^{j_n} \simeq(n/ \log n)^{1/(2t+1)}$, that
\begin{eqnarray*}
\sup_{y \in\mathbb R} |p_n(y) - p_0(y)| = \mathrm{O}\biggl(\biggl(\frac{\log
n}{n}\biggr)^{t/(2t+1)} \biggr), \qquad   \mbox{ both a.s. and in } L^p(P).
\end{eqnarray*}
For the following central limit theorem, we denote by $\rightsquigarrow
_{\ell^{\infty}(\mathbb R)}$ convergence in law for sample-bounded
processes in the Banach space of bounded functions on $\mathbb R$, and
by $G_P$ the usual $P$-Brownian bridge (for example, Chapter 3 in
\cite{1999Dudley}).
We should emphasize that the optimal bandwidth choice $2^{-j_n} \simeq
n^{-1/(2t+1)}$ (if sup-norm loss is being considered, replace $n$ by $n/
\log n$) is admissible for every $t>0$ in the theorem below.

\begin{theorem} \label{clt0}
Assume that the density $p_0$ of $P$ is a bounded function ($t=0$) or
that $p_0 \in\mathsf C^t(\mathbb R)$ for some $t$, $0 < t \le r$. Let
$j_n$ satisfy $n/(2^{j_n}j_n) \to\infty$ and $\sqrt n 2^{-j_n(t+1)}
\to0$ as $n \to\infty$. If $F$ is the distribution function of $P$
and we set $F_n^S(s) := \int_{- \infty}^s p(y,j_n)\,\mathrm{d}y$, then
\[
\sqrt n (F_n^S-F) \rightsquigarrow_{\ell^{\infty}(\mathbb R)} G_P.
\]
\end{theorem}

\begin{pf}
Given $\varepsilon>0$, apply Proposition \ref{spldf} below with
$\lambda= \varepsilon$ so that $\|F_n^S-F_n\|_\infty=\mathrm{o}_P(1/\sqrt n)$
follows and use the fact that $\sqrt n (F_n-F)$ converges in law in
$\ell^{\infty}(\mathbb R)$ to $G_P$.
\end{pf}\vspace*{-3pt}

%s3 ###
\section{The adaptive estimation procedures}\label{sec3}\vspace*{-2pt}

In this section, we construct data-driven choices of the resolution
level $j$ and state the main adaptation results. As mentioned in the
\hyperref[sec1]{Introduction}, we will use Rademacher symmetrization for this. Generate
a Rademacher sequence $\varepsilon_i$, $i=1,\ldots,n$, independent of
the sample (that is, $\varepsilon_i$ takes values $1,-1$ with
probability $1/2$) and set, for $j < l$,\vspace*{-2pt}
\begin{eqnarray} \label{rademacher}
R(n,j) &=& 2 \Biggl\|\frac{1}{n} \sum_{i=1}^n \varepsilon_i K_{j}(X_i, \cdot) \Biggr\|
_\infty   \quad \mbox{and}\nonumber
\\[-9pt]\\[-9pt]
 T(n,j,l) &=& 2 \Biggl\|\frac{1}{n} \sum_{i=1}^n
\varepsilon_i (K_{j}-K_l)(X_i, \cdot) \Biggr\|_\infty,\nonumber\vspace*{-2pt}
\end{eqnarray}
where $K_j$ is the kernel of the wavelet projection $\pi_j$ onto $V_j$
(both for Battle--Lemari\'e and compactly supported wavelets). In both
cases, these are suprema of fixed random functions that depend only on
known quantities that can be computed in a numerically effective way.
For more details on Rademacher processes, see Section \ref{radsec}.

To construct the estimators, we first need a grid indexing the spaces
$V_j$ onto which we project~$P_n$. For $r \geq1$, $n>1$, choose
integers $j_{\min}:= j_{\min, n}$ and $j_{\max}:= j_{\max, n}$ such
that $0<j_{\min} < j_{\max}$,\vspace*{-2pt}
%
%e7 ###
\begin{equation} \label{range}
2^{j_{\min}} \simeq\biggl(\frac{n}{\log n} \biggr)^{1/(2r+1)} \quad  \mbox{and}\quad
2^{j_{\max}} \simeq\frac{n}{(\log n)^2} ,\vspace*{-2pt}
\end{equation}
and set\vspace*{-2pt}
\[
\mathcal J := \mathcal J_n = [j_{\min}, j_{\max}] \cap\mathbb N.\vspace*{-2pt}
\]
Note that the number of elements in this grid is of order $\log n$. We
will consider two preliminary estimators, $\bar j_n$ and $\tilde j_n$,
of the resolution level (of course, only one is needed, but we offer a
choice between two, as discussed below). Let $p_n(j)$ be as in (\ref
{est0}) or (\ref{est}). First, we set\vspace*{-2pt}
\begin{eqnarray}\label{htf2}
\bar j_n = \min \Biggl\{ j \in\mathcal J&\dvtx& \| p_n(j) - p_n(l) \|_\infty\nonumber
\\[-10pt]\\[-10pt]
&&\quad\ \le
T(n,j,l) + 7 \|\Phi\|_2 \|p_n(j_{\max})\|^{1/2}_\infty\sqrt{\frac{2^l l}{n}},
\forall l>j, l\in{\mathcal J} \Biggr\},\nonumber
\end{eqnarray}

\noindent
where the function $\Phi$ is as in (\ref{major}), and we discuss an
explicit way to construct $\Phi$ in Remark \ref{const1} below. If the
minimum does not exist, then we set $\bar j_n$ equal to $j_{\max}$. An
alternative estimator of the resolution level is
\begin{eqnarray}\label{htf3}
\tilde j_n &=& \min \Biggl\{ j \in\mathcal J: \| p_n(j) - p_n(l) \|_\infty
\le\bigl(B(\phi)+1\bigr) R(n,l)\nonumber
\\[-8pt]\\[-8pt]
&&{}\qquad\hspace*{122pt} + 7 \|\Phi\|_2 \|p_n(j_{\max})\|^{1/2}_\infty
\sqrt{\frac{2^l l}{n}},
\forall l>j, l\in{\mathcal J} \Biggr\},\nonumber\qquad
\end{eqnarray}
where $B(\phi)$ is a bound, uniform in $j$, for the operator norm in
$L^\infty(\mathbb R)$ of the projection $\pi_j$; see Remark \ref
{const2} below. Again, if the minimum does not exist, we set $\tilde
j_n$ equal to $j_{\max}$.

Before we state the main result, we briefly discuss these procedures.
The data-driven resolution level $\tilde j_n$ in (\ref{htf3}) is based
on tests that use Rademacher-type analogs of the usual thresholds in
Lepski's method: starting with $j_{\min}$, the main contribution to $\|
p_n(j)-p_n(l)\|_\infty$ is the bias $\|Ep_n(j)-p_0\|_\infty$. The
procedure should stop when the `variance term' $\|p_n(l)-Ep_n(l)\|
_\infty$ starts to dominate. Since this is an unknown quantity and
since we know no good non-random upper bound for it, we estimate it by
the supremum of the associated Rademacher process, that is, by
$R(n,l)$. The constant $B(\phi)$ is necessary in order to correct for
the lack of monotonicity of the $R(n,l)$'s in the resolution level $l$.

The estimator $\bar j_n$ in (\ref{htf2}) is somewhat more refined: it
attempts to take advantage of the fact that in the `small bias' domain,
and using the results from Section \ref{radsec},
\[
\|p_n(j)-p_n(l)\|_\infty= \Biggl\|\frac{1}{n} \sum_{i=1}^n (K_{j}-K_l)(X_i,
\cdot) \Biggr\|_\infty
\]
should not exceed its Rademacher symmetrization
\[
T(n,j,l)=2\Biggl\|\frac{1}{n} \sum_{i=1}^n \varepsilon_i (K_{j}-K_l)(X_i,
\cdot) \Biggr\|_\infty.
\]
We now state the main result, whose proof is deferred to the next
section. As usual, we say that a wavelet basis is $s$\emph{-regular},
$s \in\mathbb{N}\cup\{0\}$, if either the scaling function $\phi$ has
$s$ weak derivatives contained in $L^p(\mathbb R)$ for some $p \ge1$
or if the mother wavelet $\psi$ satisfies $\int x^\alpha\psi(x)\,\mathrm{d}x =0$
for $\alpha= 0,\ldots,s$. Note that any compactly supported element of
$\mathsf C^s(\mathbb R), 0<s\leq 1,$ is of bounded $(1/s)$-variation so that
the $p$-variation condition in the following theorem is satisfied, for
example, for all Daubechies wavelets. The estimators below achieve the
optimal rate of convergence for estimating $p_0$ in sup-norm loss in
the minimax sense (over H\"{o}lder balls); see, for example, \cite{1999Korostelev} for optimality of these rates.

\begin{theorem} \label{cdf}
Let $X_1,\ldots,X_n$ be i.i.d. on $\mathbb{R}$ with common law $P$ that
possesses a uniformly continuous density $p_0$. Let $p_n(j):=p_n(y, j)$
be as in \textup{(\ref{est0})}, where $\phi$ is either compactly supported, of
bounded $p$-variation ($p\geq1$) and ($r-1$)-regular or $\phi= \phi
_r$ equals a Battle--Lemari\'{e} wavelet. Let the sequence $\{\hat j_n\}
_{n \in\mathbb N}$ be either $\{\bar j_n\}_{n \in\mathbb N}$ or $\{
\tilde j_n\}_{n \in\mathbb N}$ and let $F_n(\hat j_n)(t) = \int
_{-\infty}^t p_n (y,\hat j_n)\,\mathrm{d}y$. Then
%
%e8 ###
\begin{equation}\label{clt3}
\sqrt n \bigl( F_n(\hat j_n) - F \bigr) \rightsquigarrow_{\ell^\infty(\mathbb
{R})} G_P,
\end{equation}
the convergence being uniform over the set of all probability measures
$P$ on $\mathbb R$ with densities $p_0$ bounded by a fixed constant, in
any distance that metrizes convergence in law. Furthermore, if $C$ is
any precompact subset of $\mathsf C(\mathbb R)$, then
%
%e9 ###
\begin{equation}\label{sup0}
\sup_{p_0 \in C}E \sup_{y \in\mathbb R} |p_n(y, \hat j_n) - p_0(y)| = \mathrm{o}(1).
\end{equation}
If, in addition, $p_0 \in\mathsf C^t(\mathbb{R})$ for some $0<t \leq
r$, then we also have
%
%e10 ###
\begin{equation}\label{sup}
\sup_{p_0\dvtx\|p_0\|_{t, \infty} \leq D}E\sup_{y \in\mathbb R} |p_n(y,
\hat j_n) - p_0(y)| = \mathrm{O}\biggl(\biggl(\frac{\log n}{n}\biggr)^{t/(2t+1)} \biggr).
\end{equation}
\end{theorem}

\begin{remark}[(Relaxing the uniform continuity assumption)] \label{ap}
The assumption of uniform continuity of the density of $F$ can be
relaxed by modifying the definition of $\bar j_n$ (or $\tilde j_n$)
along the lines of~\cite{2009aGin}. The idea is to constrain
all candidate estimators to lie in a ball of size $\mathrm{o}(1/\sqrt n)$ around
the empirical distribution function $F_n$ so that (\ref{clt3}) holds
automatically. Formally, this can be done by adding the requirement
\[
\sup_{t \in\mathbb R} \biggl|\int_{- \infty}^t p_n(y,j)\,\mathrm{d}y - F_n(t) \biggr| \le
\frac{1}{\sqrt n \log n}
\]
to each test in (\ref{htf2}) or (\ref{htf3}). If this requirement does
not even hold for $j_{\mathrm{max}}$, then it can be seen as evidence that $F$
has no density and one just uses $F_n$ as the estimator so as to obtain
at least the functional CLT. If $F$ has a bounded density, then one can
use the exponential bound in Proposition \ref{spldf} in the proof to
control rejection probabilities of these test in the `small bias'
domain $\hat j_n > j^*$ and Theorem \ref{cdf} can then still be proven
for this procedure without \textit{any} assumptions on~$F$. See Theorem
2 in \cite{2009aGin} for more details on this procedure and
its proof.
\end{remark}

\begin{remark}[(The constant $\bolds{\|\Phi\|}_\mathbf{2}$)]\label{const1}
Once the wavelet $\phi$ have been chosen, $\hat j_n$ is purely
data-driven since the function $\Phi$ depends only on $\phi$. For the
Haar basis ($\phi=I_{[0,1)}$), we can take $\Phi=\phi$ because, in this case,
$K(x,y)\le I_{[0,1)}(|x-y|)$ so that $\|\Phi\|_2=1$. A general way to
obtain majorizing kernels $\Phi$ is described in Section 8.6 of \cite{1998Har}.
For Battle--Lemari\'{e} wavelets, the spline representation of
the projection kernel is again useful for estimating $\|\Phi\|_2$. See
\cite{1993Huang} for explicit computations.
\end{remark}

\begin{remark}[(The constant $\bolds{B(\phi)}$)]\label{const2}
To construct $\tilde j_n$, one requires knowledge of the constant
$B(\phi)$ that bounds the operator norm $\|\pi_j\|'_\infty$ of $\pi_j$,
viewed as an operator $L^\infty(\mathbb R)$. A simple way of obtaining
a bound is as follows: for any $f\in L^{\infty}(\mathbb R)$, we have,
by (\ref{major}),
\[
|\pi_j(f)(x)|=\biggl|\int K_j(x,y)f(y)\,\mathrm{d}y\biggr|\le\|\Phi\|_1\|f\|_\infty,
\]
that is, $\|\pi_j\|'_\infty\le\|\Phi\|_1.$ In combination with the
previous remark, one readily obtains possible values for $B(\phi)$. For
instance, for the Haar wavelet, $B(\phi)\le1$. For spline wavelets,
other methods are available. For example, for Battle--Lemari\'{e}
wavelets arising from linear $B$-splines, $\|\pi_j\|'_\infty$ is
bounded by $3$, and \cite{2001Shadrin}, page 135, conjectures the bound
$2r-1$ for general order $r$. See  \cite{1993DeVore}, Chapter
13.4, \cite{2001Shadrin} and references therein for more information.
\end{remark}

We also note that -- as the results in Section \ref{radsec}, in
particular Proposition \ref{random0}, show -- all of our proofs go
through if one replaces $R(n,j)$, $T(n,j,l)$ by their respective
Rademacher expectations $E^\varepsilon R(n,j)$, $E^\varepsilon
T(n,j,l)$ in the definitions of $\tilde j_n$, $\bar j_n$.

%s3.1 ###
\subsection{Estimating suprema of empirical processes}

Talagrand's \cite{1996Talagrand} exponential inequality for empirical processes (see
also \cite{2001Ledoux}), which is a uniform Prohorov-type inequality, is
not specific about constants. Constants in its Bernstein-type version
have been specified by several authors
\cite{2000Massart,2003Bousquet,2005Klein}. Let $X_i$ be the coordinates of the product
probability space $(S,\mathcal{S},P)^{\mathbb N}$, where $P$ is any
probability measure on $(S, \mathcal S)$ and let $\mathcal F$ be a
countable class of measurable functions on $S$ that take values in
$[-1/2,1/2]$ or, if $\mathcal F$ is $P$-centered, in $[-1,1]$. Let
$\sigma\le1/2$ and $V$ be any two numbers satisfying
%
%e11 ###
\begin{equation}\label{v}
\sigma^2 \geq\|Pf^2\|_\mathcal{F},\qquad V \geq n\sigma^2+2E\Biggl\|\sum
_{i=1}^n\bigl(f(X_i)-Pf\bigr)\Biggr\|_\mathcal{F},
\end{equation}
in which case $V$ is also an upper bound for $E\|\sum(f(X_i)-Pf)^2\|
_\mathcal{F}$ \cite{2005Klein}.
Then, noting that $\sup_{f\in\mathcal{F}\cup(-\mathcal{F})}\sum
_{i=1}^nf(X_i)=\sup_\mathcal{F}|\sum_{i=1}^nf(X_i)|$, Bousquet's
\cite{2003Bousquet}
version of Talagrand's inequality is as follows: for every $t >0$,
%
%e12 ###
\begin{equation}\label{bous}
\Pr\Biggl\{\Biggl\|\sum_{i=1}^n\bigl(f(X_i)-Pf\bigr)\Biggr\|_\mathcal{F}\ge E\Biggl\|\sum
_{i=1}^n\bigl(f(X_i)-Pf\bigr)\Biggr\|_\mathcal{F}+t\Biggr\}\le\exp\biggl(-\frac{t^2}{2V+(2/3)t}\biggr).
\end{equation}
In the other direction, the Klein and Rio \cite{2005Klein} result is that for
every $t>0$,
%
%e13 ###
\begin{equation}\label{kr}
\Pr\Biggl\{\Biggl\|\sum_{i=1}^n\bigl(f(X_i)-Pf\bigr)\Biggr\|_\mathcal{F}\le E\Biggl\|\sum
_{i=1}^n\bigl(f(X_i)-Pf\bigr)\Biggr\|_\mathcal{F}-t\Biggr\}\le\exp\biggl(-\frac{t^2}{2V+2t}\biggr).
\end{equation}

These inequalities can be applied in conjunction with an estimate of
the expected value obtained via empirical process methods. Here, we
describe one such result for VC-type classes, that is, for $\mathcal F$
satisfying the uniform metric entropy condition
%
%e14 ###
\begin{equation}\label{vcc}
\sup_QN(\mathcal{F}, L^2(Q),\tau)\le\biggl(\frac{A}{\tau}\biggr)^v,\qquad
0<\tau\le1\,
(A\ge e, v\ge2),
\end{equation}
with the supremum extending over all Borel probability measures on
$(S,\mathcal{S})$. We denote here by $N(\mathcal G, L^2(Q), \tau)$ the
usual covering numbers of a class $\mathcal G$ of functions by balls of
radius less than or equal to $\tau$ in $L^2(Q)$-distance. One then has,
for every $n$,
%
%e15 ###
\begin{equation}\label{expect}
E\Biggl\|\sum_{i=1}^n\bigl(f(X_i)-Pf\bigr)\Biggr\|_\mathcal{F}\le2\Biggl[15\sqrt{2vn\sigma^2\log\frac
{5A}{\sigma}}+1350v\log\frac{5A}{\sigma}\Biggr];
\end{equation}
see Proposition 3 in \cite{2009aGin} with a change obtained by
using $V$ as in (\ref{v}) instead of an earlier bound due to Talagrand
for $E\|\sum(f(X_i)-Pf)^2\|_\mathcal{F}$. Inequalities of this type also
have some historical precedents (\cite{1994Talagrand,2000Einmahl,2001Gin,2006Gin}
among others). The constants on the right-hand side of (\ref{expect})
may be far from the best possible, but we prefer them over unspecified
`universal' constants.

As is the case of Bernstein's inequality in $\mathbb R$, Talagrand's
inequality is especially useful in the Gaussian tail range and,
combining (\ref{bous}) and (\ref{expect}), one can obtain such a
`Gaussian tail' bound for the supremum of the empirical process that
depends only on $\sigma$ (similar to a bound in \cite{2001Gin}).

\begin{proposition}\label{est1}
Let $\mathcal F$ be a countable class of measurable functions that
satisfies (\ref{vcc}) and is uniformly bounded (in absolute value) by
1$/$2. Assume, further, that for some $\lambda>0$,
%
%e16 ###
\begin{equation}\label{sigmab}
n \sigma^2\ge\frac{\lambda^2v}{2}\log\frac{5A}{\sigma}.
\end{equation}
Set $c_1(\lambda) = 2[15+ 1350\lambda^{-1}]$ and let $c_2(\lambda) \geq
1+120\lambda^{-1} + 10{,}800\lambda^{-2}$. Then, if
%
%e17 ###
\begin{equation}\label{tb}
c_1(\lambda) \sqrt{2vn\sigma^2\log\frac{5A}{\sigma}}\le t\le\frac
{3}{2} c_2(\lambda) n\sigma^2,
\end{equation}
we have
%
%e18 ###
\begin{equation}\label{gaussb}
\Pr\Biggl\{\Biggl\|\sum_{i=1}^n\bigl(f(X_i)-Pf\bigr)\Biggr\|_\mathcal{F}\ge2t\Biggr\}\le\exp\biggl(-\frac
{t^2}{3c_2(\lambda) n\sigma^2}\biggr).
\end{equation}
\end{proposition}

\begin{pf}
In the light of (\ref{sigmab}), inequality (\ref{expect}) gives
\[
E\Biggl\|\sum_{i=1}^n\bigl(f(X_i)-Pf\bigr)\Biggr\|_\mathcal{F}\le c_1(\lambda)\sqrt{2vn\sigma
^2\log\frac{5A}{\sigma}}
\]
and (\ref{v}) implies that we can take $V = c_2(\lambda) n\sigma^2.$
The result now follows from (\ref{bous}), taking into account that in
the range of $t$'s,
$E\|\sum_{i=1}^n(f(X_i)-Pf)\|_\mathcal{F}\le t\le3V/2,$ (\ref{bous}) becomes
\[
\Pr\Biggl\{\Biggl\|\sum_{i=1}^n\bigl(f(X_i)-Pf\bigr)\Biggr\|_\mathcal{F}\ge2t\Biggr\}\le\exp\biggl(-\frac{t^2}{3V}\biggr).
\]
\upqed\end{pf}

The constants here may be too large for some applications, but they are
not so in situations where $\lambda$ can be taken very large, in
particular, in asymptotic considerations. (Then $c_1(\lambda) \to30$
and $c_2(\lambda) \to1$ as $\lambda\to\infty$.)

%s3.1.1 ###
\subsubsection{Estimating the size of empirical processes by
Rademacher averages}\label{radsec}\vspace*{2pt}

The constants one could obtain from Proposition \ref{est1} are not
satisfactory for the applications to adaptive estimation which we have
in mind. We now propose a remedy for this problem, inspired by a nice
idea of Koltchinskii \cite{2001Koltchinskii} and Bartlett, Boucheron and Lugosi
\cite{2002Bartlett}
which they used in other contexts, namely in risk minimization and
model selection. This consists of replacing the expectation of the
supremum of an empirical process by the supremum of the associated
Rademacher process. An inequality of this type (see \cite{2006Koltchinskii}, page 2602)
is\looseness=1\vspace*{2pt}
%
%e19 ###
\begin{equation}\label{kolt}
\Pr\Biggl\{ \Biggl\|\sum_{i=1}^n\bigl(f(X_i)-Pf\bigr)\Biggr\|_\mathcal{F}\ge2\Biggl\|\sum_{i=1}^n\varepsilon
_if(X_i)\Biggr\|_\mathcal{F}+3t\Biggr\}\le\exp\biggl(-\frac{2t^2}{3n}\biggr),\vspace*{2pt}
\end{equation}
where $\varepsilon_i$, $i \in\mathbb N$, are i.i.d. Rademacher random
variables, independent of the $X_i$'s, all defined as coordinates on a
large product probability space. Note that this bound does not take the
variance $V$ in (\ref{bous}) into account, but in the applications to
density estimation that we have in mind, $V$ is much smaller than $n$
(it is of order $n2^{-j_n}$, $j_n \to\infty$). We need a similar
inequality, with the quantity $n$ in the bound replaced by $V$, valid
over a large enough range of $t$'s.

It will be convenient to use the following well-known symmetrization
inequality (see, for example, \cite{1999Dudley}, page 343):\vspace*{2pt}
%
%e20 ###
\begin{equation} \label{radsym}
\frac{1}{2}E\Biggl\|\sum_{i=1}^n\varepsilon_if(X_i)\Biggr\|_{\mathcal{F}} - \frac{\sqrt
n}{2} \|Pf\|_\mathcal F \le E\Biggl\|\sum_{i=1}^n\bigl(f(X_i)-Pf\bigr)\Biggr\|_\mathcal{F} \le2
E\Biggl\|\sum_{i=1}^{n}\varepsilon_if(X_i)\Biggr\|_{\mathcal{F}}.\quad\vspace*{2pt}
\end{equation}
The following exponential bound is the Bernstein-type analog of (\ref
{kolt}). Denote by $E^\varepsilon$ expectation with respect to the
Rademacher variables only.

\begin{proposition}\label{random0}
Let $\mathcal F$ be a countable class of measurable functions, uniformly
bounded (in absolute value) by 1$/$2. Then, for every $t>0$,\vspace*{2pt}
%
%e21 ###
\begin{equation} \label{random}
\Pr\Biggl\{ \Biggl\|\sum_{i=1}^n\bigl(f(X_i)-Ef(X)\bigr)\Biggr\|_\mathcal{F}\ge2\Biggl\|\sum
_{i=1}^n\varepsilon_if(X_i)\Biggr\|_\mathcal{F}+3t\Biggr\}\le2\exp\biggl(-\frac{t^2}{2V'+2t}\biggr),\vspace*{2pt}
\end{equation}
as well as\vspace*{2pt}
%
%e22 ###
\begin{equation} \label{random2}
\Pr\Biggl\{ \Biggl\|\sum_{i=1}^n\bigl(f(X_i)-Ef(X)\bigr)\Biggr\|_\mathcal{F}\ge2 E^\varepsilon\Biggl\|\sum
_{i=1}^n\varepsilon_if(X_i)\Biggr\|_\mathcal{F}+3t\Biggr\}\le2\exp\biggl(-\frac{t^2}{2V'+2t}\biggr),\vspace*{2pt}
\end{equation}
where
$V'=n\sigma^2+4E\|\sum_{i=1}^n\varepsilon_if(X_i)\|_\mathcal{F}.$
\end{proposition}

\begin{pf}
We have
\begin{eqnarray*}
&&\Pr\Biggl\{ \Biggl\|\sum_{i=1}^n\bigl(f(X_i)-Pf\bigr)\Biggr\|_\mathcal{F}\ge2\Biggl\|\sum
_{i=1}^n\varepsilon_if(X_i)\Biggr\|_\mathcal{F}+3t\Biggr\}
\\
&&\quad  \le \Pr\Biggl\{ \Biggl\|\sum_{i=1}^n\bigl(f(X_i)-Pf\bigr)\Biggr\|_\mathcal{F}\ge2E\Biggl\|\sum
_{i=1}^n\varepsilon_if(X_i)\Biggr\|_\mathcal{F}+t\Biggr\}
\\
  &&{}\qquad     +\Pr\Biggl\{\Biggl\|\sum_{i=1}^n\varepsilon_if(X_i)\Biggr\|_\mathcal{F}\le E\Biggl\|\sum
_{i=1}^n\varepsilon_if(X_i)\Biggr\|_\mathcal{F} -t\Biggr\}.
\end{eqnarray*}
For the first term, combining (\ref{radsym}) with (\ref{bous}) gives
\[
\Pr\Biggl\{ \Biggl\|\sum_{i=1}^n\bigl(f(X_i)-Pf\bigr)\Biggr\|_\mathcal{F}\ge2E\Biggl\|\sum
_{i=1}^n\varepsilon_if(X_i)\Biggr\|_\mathcal{F}+t\Biggr\} \le\exp \biggl(- \frac{t^2}{2V'
+ (2/3)t} \biggr).
\]
For the second term, note that (\ref{kr}) applies to the randomized
sums $\sum_{i=1}^n \varepsilon_i f(X_i)$ as well, by just taking the
class of functions
\[
\mathcal{G}=\{g(\tau,x)=\tau f(x)\dvtx f\in\mathcal{F}\},
\]
$\tau\in\{-1,1\}$, instead of $\mathcal F$ and the probability measure
$\bar P = 2^{-1}(\delta_{-1}+\delta_1)\times P$ instead of $P$. Hence,
%
%e23 ###
\begin{equation}\label{kr'}
\Pr\Biggl\{\Biggl\|\sum_{i=1}^n\varepsilon_if(X_i)\Biggr\|_\mathcal{F}\le E\Biggl\|\sum
_{i=1}^n\varepsilon_if(X_i)\Biggr\|_\mathcal{F}-t\Biggr\}\le\exp\biggl(-\frac{t^2}{2V'+2t}\biggr)
\end{equation}
since $V' \geq n\sigma^2+2E\|\sum_{i=1}^n\varepsilon_if(X_i)\|_{\mathcal
F}$. Combining the bounds completes the proof of (\ref{random}).

It remains to prove (\ref{random2}). Let $\mathcal G$, $\bar P$ be as
above, let $Y_i = (\varepsilon_i, X_i)$ and note that $\bar P$ is the
law of~$Y_i$. By convexity,
\[
E \mathrm{e}^{-t E^\varepsilon\|\sum_{i=1}^n\varepsilon_i f(X_i)\|_\mathcal F}
\le E \mathrm{e}^{-t \|\sum_{i=1}^n \varepsilon_i f(X_i)\|_\mathcal F} = E \mathrm{e}^{-t
\|\sum_{i=1}^n g(Y_i)\|_\mathcal G}
\]
for all $t$.
The Klein and Rio \cite{2005Klein} version (\ref{kr}) of Talagrand's inequality
is, in fact, established by estimating the Laplace transform $E \mathrm{e}^{-t \|
\sum_{i=1}^ng(Y_i)\|_\mathcal G}$ and Theorem 1.2a in \cite{2005Klein} implies that
\[
E \mathrm{e}^{-t E^\varepsilon\|\sum_{i=1}^n\varepsilon_i (f(X_i)-Pf)\|
_\mathcal F} \le-t E\Biggl\|\sum_{i=1}^ng(Y_i)\Biggr\|_\mathcal G + \frac
{V}{9}(\mathrm{e}^{3t}-3t+1 )
\]
for $V \ge n \sigma^2 + 2 E \|\sum_{i=1}^n g(Y_i)\|_\mathcal G$, which,
by their proof of the implication $(a) \Rightarrow(c)$ in that
theorem, gives
\[
\Pr\Biggl\{E^\varepsilon\Biggl\|\sum_{i=1}^n\varepsilon_if(X_i)\Biggr\|_\mathcal{F}\le E\Biggl\|
\sum_{i=1}^n\varepsilon_if(X_i)\Biggr\|_\mathcal{F} -t\Biggr\} \le\exp\biggl(-\frac{t^2}{2V'+2t}\biggr).
\]
The proof of (\ref{random2}) now follows as in the previous case.
\end{pf}

For $\mathcal F$ of VC-type, the moment bound (\ref{expect}) is usually
proved as a consequence of a bound for the Rademacher process. In fact,
the proof of Proposition 3 in \cite{2009aGin} shows
that
%
%e24 ###
\begin{equation}\label{symexp}
E\Biggl\|\sum_{i=1}^n\varepsilon_if(X_i)\Biggr\|_\mathcal{F}\le15\sqrt{2vn\sigma^2\log
\frac{5A}{\sigma}}+1350v\log\frac{5A}{\sigma},
\end{equation}
where $\sigma$ is as in (\ref{v}), which we use in the following
corollary, together with the previous proposition. The constant
$c_2(\lambda)$ in the exponent below is still potentially large, but
tends to one if $\lambda\to\infty$.

\begin{corollary}\label{est2}
Let $\mathcal F$ be a countable class of measurable functions that
satisfies (\ref{vcc}) and assume it to be uniformly bounded (in
absolute value) by 1$/$2. Assume, further,
(\ref{sigmab}) for some $\lambda>0$. Then, for
$0<t\le\frac{1}{20}c_2(\lambda) n\sigma^2$
with $c_2(\lambda)$ as in Proposition \ref{est1}, we have
\begin{eqnarray*}\label{random1}
\Pr\Biggl\{\Biggl\|\sum_{i=1}^n\bigl(f(X_i)-Ef(X)\bigr)\Biggr\|_\mathcal{F}\ge2\Biggl\|\sum
_{i=1}^n\varepsilon_if(X_i)\Biggr\|_\mathcal{F}+3t\Biggr\}\le2\exp\biggl(-\frac
{t^2}{2.1c_2(\lambda) n\sigma^2}\biggr)
\end{eqnarray*}
and the same inequality holds if $\|\sum_{i=1}^n\varepsilon_if(X_i)\|
_\mathcal{F}$ is replaced by its $E^\varepsilon$ expectation.
\end{corollary}

\begin{pf}
By (\ref{sigmab}) and (\ref{symexp}), we have $V' \le c_2(\lambda) n
\sigma^2$, and the condition on $t$ together with (\ref{random}) gives
the result.
\end{pf}

%s3.2 ###
\subsection{Projections onto spline spaces and their wavelet representation} \label{splines}

In this section, we briefly review how the wavelet estimator (\ref
{est0}) for Battle--Lemari\'{e} wavelets can be represented as a spline
projection estimator (\ref{est}). We shall need the spline
representation in some proofs, while the wavelet representation will be
useful in others.

Let $T := T_j = \{t_i(j)\}_{-\infty}^\infty= 2^{-j} \mathbb Z$, $j \in
\mathbb Z$, be a bi-infinite sequence of equally spaced knots, $t_i :=
t_i(j)$. A function $S$ is a spline of order $r$, or of degree $m=
r-1$, if, on each interval $(t_i, t_{i+1})$, it is a polynomial of
degree less than or equal to $m$ (and of degree exactly $m$ on at least
one interval) and, at each breakpoint $t_i$, $S$ is at least
$(m-1)$-times differentiable. The Schoenberg space $\mathcal S_r (T) :=
\mathcal S_r (T, \mathbb R)$ is defined as the set of all splines of
order (less than or equal to) $r$ and it coincides with the space
$\mathcal S_r(T, 1, \mathbb R)$ in  \cite{1993DeVore}, page 135.
The space $\mathcal S_r (T_j)$ has a Riesz basis formed by $B$-splines
$\{N_{j,k,r}\}_{k \in\mathbb Z}$ that we now describe; see Section 4.4
in \cite{1993Schumaker} and page 138f in  \cite{1993DeVore} for more
details. Define
\[
N_{0,r}(x) = 1_{[0,1)} \ast\cdots\ast1_{[0,1)}(x),\qquad  r\mbox{-times
}      := \sum_{i=0}^{r} \frac{(-1)^i {r \choose i}
(x-i)_+^{r-1}}{(r-1)!} .
\]
For $r=2$, this is the linear $B$-spline (the usual `hat' function), for
$r=3,$ it is the quadratic and for $r=4$, it is the cubic $B$-spline. Set
$N_{k,r} (x) := N_{0,r} (x-k)$. The elements of the Riesz basis are
then given by
\[
N_{j,k, r} (x) := N_{k,r} (2^jx) = N_{0,r} (2^jx-k ).
\]
By the Curry--Schoenberg theorem, any $S \in\mathcal S_r(T_j)$ can be
uniquely represented as $S(x) = \sum_{k \in\mathbb Z} c_k
N_{j,k,r}(x).$ The orthogonal projection $\pi_j(f)$ of $f \in
L^2(\mathbb R)$ onto $\mathcal S_r(T_j) \cap L^2(\mathbb R)$ is
derived, for example, in  \cite{1993DeVore}, page 401f, where it
is shown that $\pi_j(f) = 2^{j/2} \sum_{k \in\mathbb Z} c_k
N_{j,k,r}$, with the coefficients $c_k := c_k(f)$ satisfying $(Ac)_k =
2^{j/2} \int N_{j,k,r}(x) f(x)\, \mathrm{d}x$, the matrix $A$ being given by
%
%e25 ###
\begin{equation} \label{toep}
a_{kl} = \int2^j N_{j,k,r}(x) N_{j,l,r}(x) \,\mathrm{d}x = \int N_{k,r}(x) N_{l,r}(x)\,\mathrm{d}x.
\end{equation}
The inverse $A^{-1}$ of $A$ exists (see Corollary 4.2 on page 404 in
  \cite{1993DeVore}) and if we denote its entries by $b_{kl}$ so
that $c_k = 2^{j/2} \int\sum_l b_{kl} N_{j,l,r}(x) f(x)\,\mathrm{d}x,$ then we
have
\[
\pi_j(f)(y) = 2^{j} \int\sum_{k} \sum_{l} b_{kl} N_{j,l,r}(x)
N_{j,k,r}(y) f(x)\,\mathrm{d}x = \int\kappa_j(x,y) f(x) \,\mathrm{d}x,
\]
where $\kappa_j(x,y) = 2^{j} \kappa(2^jx,2^jy)$ with
%
%e26 ###
\begin{equation} \label{prok}
\kappa(x,y) = \sum_k \sum_l b_{kl} N_{l,r} (x) N_{k,r}(y)
\end{equation}
is the spline projection kernel. Note that $\kappa$ is symmetric in its
arguments.

In fact, diagonalization of the kernel $\kappa$ of the projection
operator $\pi_j$ led
to one of the first examples of wavelets; see, for example, page 21f
and Section 2.3 in \cite{1992Meyer}, Section 5.4 in \cite{1992Daubechies} or
Section 6.1 in \cite{1998Har}. There, it is shown that there exists an
$(r-1)$-times differentiable scaling function $\phi_r$ with exponential
decay, the Battle--Lemari\'{e} wavelet of order $r$, such that
\[
\mathcal S_r(T_j) \cap L^2(\mathbb R) = V_{j,r} = \biggl\{\sum_k c_k
2^{j/2}\phi_r\bigl(2^j(\cdot)-k\bigr)\dvtx \sum_k c^2_k < \infty \biggr\}.
\]
This necessarily implies that the kernels $\kappa$ and $K=K(\phi_r)$
describe the same projections in $L^2(\mathbb R)$ and the following
simple lemma shows that these kernels are, in fact, pointwise the same.
\begin{lemma}\label{pointwise}
Let $ \{N_{k,r}\}_{k \in\mathbb Z}$ be the Riesz basis of $B$-splines
of order $r \geq1$ and let $\phi_r$ be the associated Battle--Lemari\'
{e} scaling function. If $K$ is as in (\ref{kernel0}) and $\kappa$ is
as in (\ref{prok}), then, for all $x,y \in\mathbb R$, we have
\[
K (x,y) = \kappa(x,y).
\]
\end{lemma}

\begin{pf}
If $r=1$, then $N_{0,1} = \phi_1$ since this is just the Haar basis.
So, consider $r > 1$. Since $\{\phi_r (\cdot-k) \dvtx k \in\mathbb Z \}$
is an orthonormal basis of $\mathcal S_r(\mathbb Z) \cap L^2(\mathbb
R)$ (see, for example, Theorem 1 on page 26 in \cite{1992Meyer}), it
follows that $K$ and $\kappa$ are the kernels of the same
$L^2$-projection operator and, therefore, for all $f,g \in L^2(\mathbb R),$
\[
\int\int\bigl(K (x,y) - \kappa(x,y)\bigr) f(x) g(y)\,\mathrm{d}x\,\mathrm{d}y =0.
\]
By density in $L^2(\mathbb R \times\mathbb R)$ of linear combinations
of products of elements of $L^2(\mathbb R)$, this implies that $\kappa$
and $K$ are almost everywhere equal in $\mathbb R^2$. We complete the
proof by showing that both functions are continuous on $\mathbb R^2$.
For $K$, this follows from the decomposition
\begin{eqnarray*}
|K(x,y)-K(x',y')| &\le&\sum_k |\phi_r(x-k)-\phi_r(x'-k)| |\phi_r(y-k)|
\\
 &&{}  + \sum_k |\phi_r(y-k) - \phi
_r(y'-k)| |\phi_r(x'-k)|,
\end{eqnarray*}
the uniform continuity of $\phi_r$ ($r>1$) and relation (\ref{major}).
For $\kappa$, we use the relation (\ref{kernel}) below,
\begin{eqnarray*}
|\kappa(x,y)-\kappa(x',y')| &\le&\sum_i |N_{i,r}(x)-N_{i,r}(x')| |H(y-i)|
\\
 &&{}  + \sum_i |H(y-i) - H(y'-i)| |N_{i,r}(x')|,
\end{eqnarray*}
which implies continuity of $\kappa$ on $\mathbb R^2$ since $N_{0,r}$
and $H$ are uniformly continuous (as $N_{0,r}$ is, and $\sum_i |g(|i|)|
< \infty$) and since $N_{0,r}$ has compact support.
\end{pf}

%s3.3 ###
\subsection{An exponential inequality for the uniform deviations of
the linear estimator} \label{expineq}

To control the uniform deviations of the linear estimators from their
means, one can use inequalities for the empirical process indexed by
classes of functions $\mathcal F$ contained in
%
%e27 ###
\begin{equation} \label{quai}
\mathcal K = \bigl\{2^{-j} K_j(\cdot,y)\dvtx y \in\mathbb R, j \in\mathbb N
\cup\{0\} \bigr\},
\end{equation}
together with suitable bounds on the `weak' variance $\sigma$.

If $\phi$ has compact support (and is of finite $p$-variation), it is
proved in Lemma 2 of \cite{2009bGin} that the class
$\mathcal K$ also satisfies the bound (\ref{vcc}). However, the proof
there does not apply to Battle--Lemari\'{e} wavelets. A different
proof, using the Toeplitz and band-limited structure of the spline
projection kernel, still enables us to prove that these classes of
functions are of Vapnik--Chervonenkis type.

\begin{lemma} \label{vc}
Let $\mathcal K$ be as in (\ref{quai}), where $\phi_r$ is a
Battle--Lemari\'{e} wavelet for some $r \geq1$. There then exist
finite constants $A \geq2$ and $v \geq2$ such that
\[
\sup_Q N(\mathcal K, L^2(Q), \varepsilon) \le\biggl(\frac{A}{\varepsilon} \biggr)^v
\]
for $0<\varepsilon< 1$ and where the supremum extends over all Borel
probability measures on $\mathbb R$.
\end{lemma}

\begin{pf}
In the case $r=1$, $\phi_1$ is just the Haar wavelet, in which case the
result follows from Lemma 2 of \cite{2009bGin}. Hence, we
assume that $r \geq2$.

The matrix $A$ is Toeplitz since, by a change of variables in (\ref
{toep}), $a_{kl} = a_{k+1,l+1}$ for all $k,l \in\mathbb Z$, and it is
band-limited because $N_{0,r}$ has compact support. It follows that
$A^{-1}$ is also Toeplitz and we denote its entries by $b_{kl} =
g(|k-l|)$ for some function $g$. Furthermore, it is known (for example,
Theorem 4.3 on page 404 of \cite{1993DeVore}) that the entries
of the inverse of any positive definite band-limited matrix satisfy
$|b_{kl}| \le c \lambda^{|k-l|}$ for some $0< \lambda< 1$ and $c$ finite.
Now, following  \cite{1993Huang}, we write
\[
\sum_k g(|l-k|)N_{k,r}(x) = \sum_{k} g(|l-k|) N_{k-l,r}(x-l) = \sum_k
g(|k|)N_{k,r}(x-l),
\]
so that
%
%e28 ###
\begin{equation} \label{kernel}
2^{-j} \kappa_j (\cdot, y) = \sum_{l \in\mathbb Z} N_{j,l,r}(y)
H\bigl(2^j(\cdot)-l\bigr) ,
\end{equation}
where $H(x) = \sum_{k \in\mathbb Z} g(|k|)N_{k,r}(x)$ is a function of
bounded variation. To see the last claim, note that $N_{0,r}$ is of
bounded variation and hence $\|N_{k,r}\|_{\mathrm{TV}} = \|N_{0,r}\|_{\mathrm{TV}}$
(where $\|\cdot\|_{\mathrm{TV}}$ denotes the usual total variation norm) so that
$\|H\|_{\mathrm{TV}} \le\|N_{0,r}\|_{\mathrm{TV}} \times\sum_{k \in\mathbb Z} |g(|k|)|
< \infty$ because $\sum_k |b_{l, l-k}| \le\sum_k c \lambda^{|k|} <
\infty$. The last fact implies that
\[
\mathcal H = \bigl\{H\bigl(2^j(\cdot)-l\bigr)\dvtx l \in\mathbb Z, j \in\mathbb N \cup\{
0\} \bigr\}
\]
satisfies, for finite constants $B>1$ and $w \geq1,$
\[
\sup_Q N(\mathcal H, L^2(Q), \varepsilon) \le\biggl( \frac{B\|H\|_\mathrm{TV}
}{\varepsilon} \biggr)^w  \qquad \mbox{for }    0<\varepsilon< \|H\|_\infty,
\]
as proved in \cite{1987Nolan}. Since $N_{j,0,r}$ is zero if $y$
is not contained in $[0, 2^{-j}r]$, the sum in~(\ref{kernel}), for
fixed $y$ and $j$, extends over only those $l$'s such that $2^{j}y -r
\le l < 2^jy$, hence it consists of at most $r$ terms. This implies
that $\mathcal K$ is contained in the set $\mathcal H_r$ of linear
combinations of at most $r$ functions from $\mathcal H$, with
coefficients bounded in absolute value by $\|N_{j,l,r}\|_\infty= \|
N_{0,r}\|_\infty< \infty$. Given $\varepsilon$, let $\varepsilon' =
\varepsilon/ (2r\max(\|H\|_\infty, \|N_{0,r}\|_\infty))$. Let $\alpha
_1,\ldots, \alpha_{n_1}$ be an $\varepsilon'$-dense subset of $[- \|
N_{0,r}\|_\infty, \|N_{0,r}\|_\infty]$ which, for $\varepsilon' < \|
N_{0,r}\|_\infty$, has cardinality $n_1 \leq3 \|N_{0,r}\|_\infty
/\varepsilon'$. Furthermore, let $h_1,\ldots,h_{n_2}$ be a subset of
$\mathcal H$ of cardinality $n_2 = N(\mathcal H, L^2(Q), \varepsilon
')$ which is $\varepsilon'$-dense in $\mathcal H$ in the
$L^2(Q)$-metric. It follows that for
$\varepsilon' < \min(\|H\|_\infty, \|N_{0,r}\|_\infty)$, every $\sum
_{l \in\mathbb Z} N_{j,l,r}(y) H(2^j(\cdot)-l)$ is at
$L^2(Q)$-distance at most $\varepsilon$ from $\sum_{l=1}^r \alpha
_{i(l)} h_{i'(l)}$ for some $1 \leq i(l) \leq n_1$ and $1 \leq i'(l)
\leq n_2$. The total number of such linear combinations is dominated by
$(n_1 n_2)^r \le(B'/\varepsilon)^{(w+1)r}$. This shows that the lemma
holds for $\varepsilon <  2r \min\{\|H\|_\infty, \|N_{0,r}\|_\infty
\}  \max\{\|H\|_\infty, \|N_{0,r}\|_\infty\} = 2r\|H\|_\infty\|
N_{0,r}\|_\infty= U$, which completes the proof by taking $A=\max
(B',U,e)$ (for $\varepsilon\in[U,A]$, one ball covers the whole set).
\end{pf}

\begin{proposition} \label{ineq}
Let $K$ be as in \textup{(\ref{kernel0})} and assume either that $\phi$ has
compact support and is of bounded $p$-variation ($p < \infty$) or that
$\phi$ is a Battle--Lemari\'{e} scaling function for some $r \geq1$.
Suppose that $P$ has a bounded density $p_0$. Given $C, T>0$, there
exist finite positive constants $C_1 = C_1(C,K, \|p_0\|_\infty)$ and
$C_2 = C_2(C,T,K, \|p_0\|_\infty)$ such that, if
\[
\frac{n}{2^j j} \geq C  \quad   \mbox{and} \quad     C_1 \sqrt{\frac{2^jj}{n}}
\leq t \leq T ,
\]
then
%
%e29 ###
\begin{equation}\label{as1}
\Pr\Bigl\{\sup_{y \in\mathbb R} |p_n(y,j)-Ep_n(y,j)| \ge t \Bigr\} \le\exp
\biggl(-C_2\frac{nt^2}{2^j}\biggr).
\end{equation}
\end{proposition}

\begin{pf} We first prove the Battle--Lemari\'{e} wavelet case. If
$r>1$, then the function $K$ is continuous (see the proof of Lemma \ref
{pointwise}) and therefore the supremum in (\ref{as1}) is over a
countable set. That this is also true for $r=1$ follows from Remark 1
in \cite{2009bGin}.
We apply Proposition \ref{est1} and Lemma \ref{vc} to the supremum of
the empirical process indexed by the classes of functions
\[
\mathcal K_j:=\{2^{-j}K_j(\cdot,y)/(2\|\Phi\|_\infty)\dvtx y\in\mathbb R\},
\]
where $\Phi$ is a function majorizing $K$ (as in (\ref{major})) so that
$\mathcal K_j$ is uniformly bounded by $1/2$. We next bound the second
moments $E (2^{-2j} K_j^2(X,y))$. We have, using (\ref{major}),
that
\begin{eqnarray} \label{secm}
\int2^{-2j} K_j^2(x,y) p_0(x)\,\mathrm{d}x &\le& \int\Phi^2\bigl(|2^{j}(x-y)|\bigr)
p_0(x)\,\mathrm{d}x \nonumber
\\[-8pt]\\[-8pt]
& \le& 2^{-j} \int\Phi^2(|u|) p_0 (y+2^{-j}u)\, \mathrm{d}u \le2^{-j} \|p_0\|
_\infty\|\Phi\|_2^2.\nonumber
\end{eqnarray}
We may hence take $\sigma= \sqrt{2^{-j} \|\Phi\|_2^2 \|p_0\|_\infty
}/(2\|\Phi\|_\infty)$ and the result is then a direct consequence of
Proposition \ref{est1}, which applies by Lemma \ref{vc}. For compactly
supported wavelets, the same proof applies, using Lemma 2 (and Remark
1) in \cite{2009bGin}.
\end{pf}

\begin{pf*}{Proof of Theorem \ref{linear}}
Using Lemma \ref{vc}, the first two claims of the theorem follow by the
same proof as in  \cite{2009bGin}, Theorem 1 and Remark~4.
For the bias term, we argue as in Theorem 8.1 in \cite{1998Har} -- using
the fact that $\phi_r$ is $(r-1)$-times differentiable -- and obtain,
for $p_0 \in\mathsf C^t(\mathbb R) ,$
%
%e30 ###
\begin{eqnarray}\label{bias0}
|Ep_n(x)-p_0(x)|&\le&2^{-jt}\|p_0\|_{t,\infty}C ,
\end{eqnarray}
where $C:=C(\Phi)= \int\Phi(|u|)|u|^t\,\mathrm{d}u$.
\end{pf*}

%s3.4 ###
\subsection{An exponential inequality for the distribution function of
the linear estimator} \label{expineq2}

The quantity of interest in this subsection is the distribution
function $F_n^S$ of the linear projection estimator $p_n$ from (\ref
{est}). More precisely, we will study the stochastic process
\[
\sqrt n \bigl(F^S_n(s)-F(s)\bigr) = \sqrt n \int_{- \infty}^s \bigl(p_n(y,j) -
p_0(y)\bigr)\, \mathrm{d}y,     \qquad  s \in\mathbb R.
\]
To prove a functional CLT for this process, it turns out that it is
easier to compare $F_n^S$ to $F_n$ rather than to $F$. With $\mathcal F
= \{1_{(-\infty, s]}\dvtx  s \in\mathbb R\}$, the decomposition
%
%e31 ###
\begin{equation} \label{dec}
(F_n^S-F_n)(s) = (P_n-P)\bigl(\pi_j (f)-f\bigr) + \int\bigl(\pi_j (p_0)-p_0\bigr) f,
   \qquad    f \in\mathcal F,
\end{equation}
will be useful, since it splits the quantity of interest into a
deterministic `bias' term and an empirical process.

\begin{lemma} \label{cltbias}
Assume that $p_0$ is a bounded function ($t=0$) or that $p_0 \in
\mathsf C^t(\mathbb R)$ for some $0<t\le r$. Let $\mathcal F = \{
1_{(-\infty, s]}\dvtx  s \in\mathbb R\}$. We then have
%
%e32 ###
\begin{equation}
\biggl|\int_\mathbb R\bigl(\pi_j(p_0)-p_0\bigr)f\biggr| \le C 2^{-j(t+1)}
\end{equation}
for some constant $C$ depending only on $r$ and $\|p_0\|_{t, \infty}$.
\end{lemma}

\begin{pf}
Let $\psi:= \psi_r$ be the mother wavelet associated with $\phi_r$.
Since the wavelet series of $p_0 \in L^1(\mathbb R)$ converges in $
L^1(\mathbb R)$, we have
$\pi_j(p_0)-p_0 = -\sum_{l=j}^\infty\sum_k \beta_{lk}(p_0) \psi_{lk}$
in the $L^1(\mathbb R)$-sense and then, since $f = 1_{(-\infty, s]} \in
L^\infty(\mathbb R)$,
\begin{eqnarray*}
-\int_{\mathbb{R}}\bigl(\pi_j(p_0)-p_0\bigr)f = \int_{\mathbb R} \Biggl(\sum_{l=j}^\infty
\sum_k \beta_{lk}(p_0) \psi_{lk}(x) \Biggr) f(x)\,\mathrm{d}x = \sum_{l=j}^{\infty} \sum
_k \beta_{lk}(p_0) \beta_{lk}(f).
\end{eqnarray*}
The lemma now follows from an estimate for the decay of the wavelet
coefficients of $p_0$ and $f$, namely, the bounds
%
%e33 ###
\begin{equation}\label{decay}
\sup_{f \in\mathcal F} \sum_k |\beta_{lk}(f)| \le c
2^{-l/2}   \quad  \mbox{and} \quad    \sup_k |\beta_{lk}(p_0)| \le c' 2^{-l(t+1/2)}.
\end{equation}
The first bound is proved as in the proof of Lemma 3 in \cite{2009bGin}, noting that the identity before equation (37) in that
proof also holds for spline wavelets by their exponential decay
property. The second bound follows from
\begin{eqnarray*}
\sup_k |\beta_{lk}(p_0)| &\le& c'' 2^{-l/2} \|K_{l+1}(p_0)-K_l(p_0)\|
_\infty
\\
&\le& c''2^{-l/2} \bigl( \|K_l(p_0)-p_0\|_\infty+ \|K_{l+1}(p_0)-p_0\|
_\infty\bigr) \le c' 2^{-l/2}2^{-lt} ,
\end{eqnarray*}
where we used (9.35) in \cite{1998Har} for the first inequality and (\ref
{bias0}) in the last.
\end{pf}

To control the fluctuations of the stochastic term, one applies
Talagrand's inequality to the empirical process indexed by the
`shrinking' classes of functions $\{\pi_j (f)-f \dvtx f \in\mathcal F\}$.
These classes consist of differences of elements in $\mathcal F$ and in
\[
\mathcal K'_j := \biggl\{\int_{- \infty} ^t K_j(\cdot,y)\, \mathrm{d}y \dvtx t \in\mathbb
R \biggr\},
\]
and we have to show that for each $j$, this class satisfies the entropy
condition (\ref{vcc}). Again, for $\phi$ with compact support (and of
finite $p$-variation), this result was proven in Lemma 2 of \cite{2009bGin} and we now extend it to the Battle--Lemari\'{e}
wavelets considered here.

\begin{lemma} \label{vc1}
Let $\mathcal K'_j$ be as above, where $\phi_r$ is a Battle--Lemari\'
{e} wavelet for $r \geq1$. There then exist finite constants $A \geq
e$ and $v \ge2$, independent of $j$ and such that
\[
\sup_Q N(\mathcal K_j', L^2(Q), \varepsilon) \le\biggl(\frac{A}{\varepsilon}
\biggr)^v,\qquad 0<\varepsilon<1,
\]
where the supremum extends over all Borel probability measures on
$\mathbb R$.
\end{lemma}

\begin{pf}
In analogy to the proof of Lemma \ref{vc}, one can write
\[
\int_{- \infty} ^t K_j(\cdot,y)\, \mathrm{d}y = \sum_{l \in\mathbb Z} \int_{-
\infty}^t 2^j N_{j,l,r}(y)\,\mathrm{d}y \,H\bigl(2^j(\cdot)-l\bigr)
\]
since the series (\ref{kernel}) converges absolutely (in view of
\[
\sum_l |H(2^jx-l)| \le\sum_k |g(|k|)| \sum_l N_{k,r} (2^jx-l) \le r \|
N_{0,r}\|_\infty\sum_k |g(|k|)| < \infty).
\]
Recall that $N_{j,l,r}$ is supported in the interval
$[2^{-j}l,2^{-j}(r+l)]$. Hence, if $l>2^jt$, then the last integral is
zero. For $l \leq2^jt-r$, the integral equals the constant $c = \int
_\mathbb R N_{0,r}(y)\, \mathrm{d}y$ and for $l \in[2^jt-r, 2^jt]$, the integral
$c_{j,l,r}$ is bounded by $c$, so this sum, in fact, equals
\[
c \sum_{l \leq2^jt -r} H\bigl(2^j(\cdot)-l\bigr) + \sum_{2^jt-r < l < 2^jt}
c_{j,l,r} H\bigl(2^j(\cdot)-l\bigr).
\]
The second sum is contained in the set $\mathcal H_r$ from the proof of
Lemma \ref{vc}, which satisfies the required entropy bound independent
of $j$. For the first sum, decompose $H$ into its positive and negative
parts, so that the two resulting collections of functions are linearly
ordered (in $t$) by inclusion and are hence a VC-subgraph of index 1;
see Theorems 4.2.6 and 4.8.1 in \cite{1999Dudley}. Moreover, we can take
the envelope $r \|N_{0,r}\|_\infty\sum_k |g(|k|)|$ independent of $j$.
Combining entropy bounds, this proves the lemma.
\end{pf}

Combining these observations, one can prove the following inequality,
which parallels Theorem 1 of \cite{2009aGin} for the
classical kernel density estimator, and Lemma 4 of \cite{2009aGin}
 for the wavelet density estimator (with $\phi$ compactly supported).

\begin{proposition} \label{spldf}
Let $F_n(s)=\int_{-\infty}^s\mathrm{d}P_n$ and $F_n^S(s):=F_n^S(s,j) = \int_{-
\infty} ^s p_n(y,j)\,\mathrm{d}y$, where $p_n$ is as in \textup{(\ref{est})}. Assume that
the density $p_0$ of $P$ is a bounded function ($t=0$) or that $p_0 \in
\mathsf C^t(\mathbb R)$ for some~$t$, $0 < t \le r$. Let $j \in\mathbb
Z$ satisfy $2^{-j} \geq d (\log n/n)$ for some $0<d<\infty$. There then
exist finite positive constants $L := L(\|p_0\|_\infty, K,d)$, $\Lambda
_0 := \Lambda_0 (\|p_0\|_{t, \infty}, K,d)$ such that for all $n \in
\mathbb N$ and $\lambda\geq\Lambda_0 \max(\sqrt{j2^{-j} }, \sqrt n
2^{-j(t+1)})$, we have
\begin{eqnarray*}
\Pr\bigl(\sqrt n \|F_n^{S} - F_n \|_\infty> \lambda\bigr) \leq L \exp \biggl\{-
\frac{\min(2^j\lambda^2, \sqrt n \lambda)}{L} \biggr\}.
\end{eqnarray*}
\end{proposition}

\begin{pf} Given the preceding lemmas, the proposition follows from
Talagrand's inequality applied to the class $\{\pi_j(1_{(-\infty
,x]})-1_{(-\infty,x]}\}$ in the same way as in the proof of Lemma 4 in
\cite{2009bGin}, so we omit it.
\end{pf}

%s3.5 ###
\subsection[Proof of Theorem 3]{Proof of Theorem \protect\ref{cdf}}

We can now prove the main result, Theorem \ref{cdf}. We will prove it
only for Battle--Lemari\'{e} wavelets. For compactly supported
wavelets, the proof is exactly the same, replacing the results from
steps~(I) and (II) below and from Sections \ref{expineq} and \ref{expineq2}
for spline wavelets by the corresponding ones for compactly supported
wavelets obtained in \cite{2009bGin}. Also, uniformity in
$p_0$ -- which is proved by controlling the respective constants -- is
left implicit in the derivations. We start with some preliminary observations.

(I) Since, uniformly in $j \in\mathcal J$, we have $n/(2^j j)
> c \log n$ for some $c>0$ independent of $n$, we have from Theorem~\ref{linear} that
%
%e34 ###
\begin{equation}\label{splvar}
E\|p_n(j) - Ep_n(j)\|_\infty^p \leq D^p \biggl({\frac{2^j j}{n}}\biggr)^{p/2} :=
D^p \sigma^p (j,n)
\end{equation}
for every $j \in\mathcal J$, $1 \le p <\infty$ and some $0<D<\infty$
depending only on $\|p_0\|_\infty$ and $\Phi$.

For the bias, we recall from (\ref{bias0}) that for $0 < t \le r$,
%
%e35 ###
\begin{equation}\label{splbias}
|E p_n(y,j) - p_0(y)| \leq2^{-jt} \|p_0\|_{t,\infty} C(\Phi) := B(j, p_0).
\end{equation}
If the density $p_0$ is only uniformly continuous, then one still has
from (\ref{major}) and integrability of $\Phi$ that, uniformly in $y
\in\mathbb R$,
%
%e36 ###
\begin{equation}\label{nonuniform}
|E p_n(y,j) - p_0(y)| \le\biggl|\int|\Phi(|u|)||p_0(y-2^{-j}u) - p_0(y)|\,\mathrm{d}u
\biggr| := B(j, p_0) = \mathrm{o}(1).
\end{equation}

(II) Define $\tilde M := \tilde M_n = C \|p_n(j_{\max})\|_\infty
$ and set $C= 49 \|\Phi\|^2_2$. Also, define $M = C \|p_0\|_\infty$ for
the same $C$. We need to control the probability that $\tilde M > 1.01
M$ or $\tilde M < 0.99 M$ if $p_0$ is uniformly continuous. For some
$0<L<\infty$ and $n$ large enough, we have
\begin{eqnarray*}
&&\Pr(|\tilde M - M | > 0.01 C\|p_0\|_\infty)
\\
&&\quad= \Pr\bigl(|\|p_n(j_{\max})\|_\infty- \|p_0\|_\infty| > 0.01 \|p_0\|
_\infty\bigr)
\\
&&\quad \leq \Pr\bigl(\|p_n(j_{\max}) - p_0\|_\infty> 0.01 \|p_0\|_\infty\bigr)
\\
&&\quad \leq \Pr\bigl(\|p_n(j_{\max}) - Ep_n(j_{\max})\|_\infty> 0.01 \|p_0\|
_\infty   - B(j_{\max}, p_0) \bigr)
\\
&&\quad \le \Pr\bigl(\|p_n(j_{\max}) - Ep_n(j_{\max})\|_\infty> 0.009 \|p_0\|
_\infty\bigr)
\\
&&\quad \le \exp\biggl\{-\frac{(\log n)^2}{L} \biggr\},
\end{eqnarray*}
by Proposition \ref{ineq} and step (I). Furthermore, there exists a
constant $L'$ such that $E \tilde M \le L'$ for every $n$, in view of
\[
E\|p_n(j_{\max})\|_\infty\le E\|p_n(j_{\max}) - Ep_n(j_{\max})\|
_\infty+ \|Ep_n(j_{\max})\|_\infty\le c + \|\Phi\|_1\|p_0\|_\infty,
\]
where we have used (\ref{major}) and (\ref{splvar}).

(III) We need some observations on the Rademacher processes used in the
definition of $\hat j_n$. First, for the symmetrized empirical measure
$\tilde P_n = 2 n^{-1} \sum_{i=1}^n \varepsilon_i \delta_{X_i}$, we have
%
%e37 ###
\begin{equation} \label{monot}
R(n,j) = \|\pi_j(\tilde P_n)\|_\infty= \|\pi_j (\pi_l(\tilde P_n))\|
_\infty\le\|\pi_j\|'_\infty R(n,l) \le B(\phi) R(n,l)
\end{equation}
for every $l>j$. Here, $\|\pi_j\|'_\infty$ is the operator norm in
$L^\infty(\mathbb R)$ of the projection $\pi_j$, which admits bounds
$B(\phi)$ independent of $j$. (Clearly, $\pi_j$ acts on finite signed
measures $\mu$ by duality, taking values in $L^\infty(\mathbb R)$ since
$|\pi_j(\mu)| =|\int K_j(\cdot,y) \,\mathrm{d}\mu(y) | \le2^j\|\Phi\|_\infty|\mu
|(\mathbb R)$.) See Remark \ref{const2} for details on how to obtain
$B(\phi)$. Furthermore, for $j < l$,
%
%e38 ###
\begin{equation} \label{triang}
T(n,j,l) \le R(n,j) + R(n,l) \le\bigl(1+B(\phi)\bigr)R(n,l)
\end{equation}
and the same inequality holds for the Rademacher expectations of
$T(n,j,l)$. We also record the following bound for the (full)
expectation of $R(n,l)$, $l \in\mathcal J$: using inequality (\ref
{symexp}) and the variance computation (\ref{secm}), we have that there
exists a constant $L$ depending only on $\|p_0\|_\infty$ and $\Phi$
such that, for every $l \in\mathcal J$, $ER(n,l) \le L \sqrt{2^l l/n}.$

\begin{pf*}{Proof of (\ref{clt3})} Let $\mathcal F = \{1_{(-\infty, s]}\dvtx s
\in\mathbb R\}$ and let $f \in\mathcal F$. We have

\[
\sqrt n \int\bigl(p_n(\hat j_n) - p_0\bigr) f = \sqrt n \int\bigl(p_n (j_{\max} ) -
p_0\bigr) f + \sqrt n \int\bigl(p_n(\hat j_n) - p_n(j_{\max})\bigr) f.
\]
The first term satisfies the CLT from Theorem \ref{clt0} for the linear
estimator with $j_n = j_{\max}$. We now show that the second term
converges to zero in probability. First, observe that
\begin{eqnarray*}
p_n(\hat j_n)(y) - p_n(j_{\max})(y) = P_n\bigl(K_{\hat j_n} (\cdot
,y)-K_{j_{\max}}(\cdot,y)\bigr) = - \sum_{l=\hat j_n}^{j_{\max}-1} \sum_k
\hat\beta_{lk} \psi_{lk}(y),
\end{eqnarray*}
with convergence in $L^1(\mathbb R)$. Next, we have, by (9.35) in \cite{1998Har}, for all $l \in[\hat j_n, j_{\max}-1]$ and all $k$, by the
definition of $\hat j_n$, that for some $0 < D' < \infty,$\vspace*{2pt}
\begin{eqnarray*}
(1/D') 2^{l/2}|\hat\beta_{lk}| &\le& \sup_{y \in\mathbb R}
|P_n(K_{l+1} (\cdot,y))- P_n(K_{l}(\cdot,y)) | = \|p_n(l+1) - p_n(l)\|
_\infty
\\[2pt]
&\le& \|p_n(l+1)-p_n(\hat j_n)\|_\infty+ \|p_n(l)-p_n(\hat j_n)\|
_\infty
\\[2pt]
&\le& \bigl(1+B(\phi)\bigr) \bigl(R(n,l+1)+ R(n,l)\bigr) + 3 \sqrt{\tilde M 2^l l/n},\vspace*{2pt}
\end{eqnarray*}
in the case $\hat j_n = \bar j_n$, also using the inequality $T(n,\bar
j_n,l) \le(1+B(\phi)) R(n,l)$ for $l \ge\bar j_n$; see (\ref
{triang}). Consequently, uniformly in $f \in\mathcal F$,\vspace*{2pt}
\begin{eqnarray*}
&& E\biggl| \int\bigl(p_n(\hat j_n) - p_n(j_{\max})\bigr) f \biggr|
\\[2pt]
&&\quad= E\Biggl| \sum_{l=\hat
j_n}^{j_{\max}-1} \sum_k \hat\beta_{lk} \int\psi_{lk}(y) f(y)\,\mathrm{d}y \Biggr|
\\[2pt]
&&\quad\le E \sum_{l=j_{\min}}^{j_{\max}-1} D' 2^{-l/2} \bigl(\bigl(B(\phi
)+1\bigr)\bigl(R(n,l+1)+R(n,l)\bigr) + 3\sqrt{\tilde M 2^l l/n} \bigr) \sum_k |\beta
_{lk}(f)|
 \\[2pt]
&&\quad \le\biggl( \frac{D''}{\sqrt n} \biggr) \sum_{l=j_{\min}}^{j_{\max}-1}
2^{-l/2}\sqrt l = \mathrm{o} \biggl(\frac{1}{\sqrt n} \biggr),\vspace*{2pt}
\end{eqnarray*}
using the moment bounds in (II) and (III), $\hat j_n \geq j_{\min} \to\infty
$ as $n \to\infty$ (by definition of $\mathcal J$) and the fact that
$\sup_{f \in\mathcal F} \sum_k |\beta_{lk}(f)| \le c2^{-l/2}$ by (\ref
{decay}) for some constant $c$.
\end{pf*}

\begin{pf*}{Proof of (\ref{sup0}) and (\ref{sup})} The proof of the case
$t=0$ follows from a simple modification of the arguments below as in
Theorem 2 of \cite{2009aGin}, so we omit it. (In this case,
one defines $j^*$ as $j_{\max}$ if $t=0$ so that only the case $\hat
j_n \le j^*$ has to be considered.) For $t>0$, define $j^*:=j(p_0)$ by
the balance equation\vspace*{2pt}
%
%e39 ###
\begin{equation} \label{star}
j^*=\min\bigl\{j\in\mathcal{J}\dvtx B(j,p_0) \le\sqrt{2 \log2} \|p_0\|
^{1/2}_\infty\|\Phi\|_2 \sigma(j,n) \bigr\}.\vspace*{2pt}
\end{equation}
Using the results from (I), it is easily verified that $2^{j^*} \simeq
(n/\log n))^{1/(2t+1)}$ if $p_0 \in\mathsf C^t(\mathbb R)$ for
some $0<t\leq r$ and that\vspace*{2pt}
\[
\sigma(j^*,n) = \mathrm{O}\biggl(\biggl(\frac{\log n}{n}\biggr)^{t/(2t+1)}\biggr)\vspace*{2pt}
\]
is the rate of convergence required in (\ref{sup}).

We will consider the cases $\{\hat j_n \leq j^* \}$ and $\{\hat j_n >
j^* \}$ separately. First, if $\hat j_n$ is $\bar j_n$, then we have,
by the definition of $\bar j_n$, (\ref{triang}), the definitions of $M$
and $j^*$, (\ref{splvar}) and the moment bound in (III),
\begin{eqnarray} \label{lowvar}
&& E \|p_n(\bar j_n) - p_0 \|_\infty I_{\{\bar j_n \le j^* \} \cap\{
\tilde M \leq1.01M\}} \nonumber
\\
&&\quad \leq E \bigl( \|p_n (\bar j_n) - p_n (j^*) \|_\infty+ E\| p_n (j^*) -
p_0\|_\infty\bigr) I_{\{\bar j_n \le j^* \} \cap\{\tilde M \leq1.01M\}}
\nonumber
\\[-8pt]\\[-8pt]
&&\quad \le \bigl(B(\phi)+1\bigr) ER(n, j^*) + \sqrt{1.01 M} \sigma(j^*,n) + \| p_n
(j^*) - p_0\|_\infty\nonumber
\\
&&\quad \le B' \sqrt{\frac{2^{j^*}j^*}{n}} + B'' \sigma(j^*,n) = \mathrm{O}(\sigma
(j^*, n)) \nonumber.
\end{eqnarray}
If $\hat j_n$ is $\tilde j_n$, then one has the same bound (without
even using (\ref{triang})).

Also, by the results in (I) and (II),
we have
\begin{eqnarray*}
&&E \|p_n(\hat j_n) - p_0 \|_\infty I_{\{\hat j_n \le j^* \} \cap\{
\tilde M > 1.01 M\}}
\\
&&\quad \le\sum_{j \in\mathcal J: j \le j^*} E \bigl( [\|p_n(j)-Ep_n(j)\|
_\infty+ B(j, p_0) ] I_{\{\hat j_n =j\}} I_{\{\tilde M > 1.01M \}} \bigr)
\\
&&\quad \le c \log n [D \sigma(j^*,n) + B(j_{\min}, p_0) ] \cdot\sqrt{E
1_{\{\tilde M > 1.01M\}}}
\\
&&\quad = \mathrm{o}\Biggl((\log n) \sqrt{\exp\biggl\{-\frac{(\log n)^2}{L} \biggr\}}\Biggr) = \mathrm{o}(\sigma
(j^*, n)) .
\end{eqnarray*}
We now turn to $\{\hat j_n > j^* \}$. First,
\begin{eqnarray*}
&& E \| p_n (\hat j_n) - p_0 \|_\infty I_{\{\hat j_n > j^* \}\cap\{
\tilde M < 0.99 M\}}
\\
&&\quad \le\sum_{j \in\mathcal J\dvtx j > j^*} E \bigl( [\|p_n(j)-Ep_n(j)\|_\infty
+ B(j, p_0) ] I_{\{\hat j_n =j\}} I_{\{\tilde M < 0.99 M \}} \bigr)
\\
&&\quad \le c' \log n [D \sigma(j_{\max}, n) + B(j^*, p_0) ] \cdot\sqrt{E
I_{\{\tilde M < 0.99 M\}}}
\\
&&\quad = \mathrm{O}\Biggl(\sqrt{(\log n) \exp\biggl\{-\frac{(\log n)^2}{L} \biggr\}}\Biggr) = \mathrm{o}(\sigma
(j^*, n)),
\end{eqnarray*}
again by the results in (I) and (II), and, second, for any $1<p<\infty$,
$1/p+1/q=1$, using (\ref{splvar}) and the definition of $j^*$, we have\vspace*{-2pt}
\begin{eqnarray*}
&& E \| p_n (\hat j_n) - p_0 \|_\infty I_{\{\hat j_n > j^* \} \cap\{
0.99M \leq\tilde M\}}
\\[-2pt]
&&\quad \leq\sum_{j \in\mathcal{J}:j > j^*} \bigl(E \|p_n (j) - p_0 \|_\infty^p
\bigr)^{1/p}   \bigl(EI_{\{\hat j_n = j\}\cap\{0.99M \leq\tilde M\}}\bigr)^{1/q}
\\[-2pt]
&&\quad \leq \sum_{j \in\mathcal{J}: j > j^*} D' \sigma(j,n) \cdot\Pr (\{
\hat j_n = j\} \cap\{0.99M \leq\tilde M\} )^{1/q}.\vspace*{-2pt}
\end{eqnarray*}
We show below that for $n$ large enough, some constant $c$, some
$\delta>0$ and some $q>1$,\vspace*{-2pt}
%
%e40 ###
\begin{equation} \label{b0}
\Pr(\{\hat j_n = j\} \cap\{0.99M \leq\tilde M\}) \le c
2^{-j(q/2+\delta)},\vspace*{-2pt}
\end{equation}
which gives the bound\vspace*{-2pt}
\[
\sum_{j \in\mathcal{J}\dvtx  j > j^*} D'' \sigma(j,n) \cdot2^{-j/2-j\delta
/q} = \mathrm{O} \biggl(\frac{1}{\sqrt n} \biggr) = \mathrm{o}(\sigma(j^*,n)),\vspace*{-2pt}
\]
completing the proof, modulo verification of (\ref{b0}).

To verify (\ref{b0}), we split the proof into two cases. Pick any $j
\in\mathcal{J}$ such that $j > j^*$ and denote by $j^-$ the previous
element in the grid (that is, $j^-= j-1$).

\textit{Case I}: $\hat j_n = \bar j_n$.
We have
\begin{eqnarray*}
&&\Pr(\{\bar j_n= j\} \cap\{0.99M \leq\tilde M\})
\\
&&\quad \leq \sum_{l \in
\mathcal{J}\dvtx l \ge j} \Pr\bigl( \| p_n(j^-) - p_n(l) \|_\infty>
T(n,j^-,l)
  + \sqrt{0.99 M} \sigma(l,n)\bigr).
\end{eqnarray*}
We first observe that
\begin{eqnarray} \label{triang2}
&&\| p_n(j^-) - p_n(l) \|_\infty\nonumber
\\[-8pt]\\[-8pt]
&&\quad \leq\| p_n(j^-) - p_n(l) - Ep_n(j^-) +
Ep_n(l) \|_\infty+ B(j^-, p_0) + B(l, p_0),\nonumber
\end{eqnarray}
where, setting $\sqrt{2 \log2}\|p_0\|^{1/2}_\infty\|\Phi\|_2 =:
U(p_0,\Phi)$,
\[
B(j^-, p_0) + B(l, p_0) \leq2B(j^*, p_0) \leq2 U(p_0,\Phi) \sigma
(j^*,n) \leq2 U(p_0, \Phi) \sigma(l,n),
\]
by definition of $j^*$ and since $l>j^- \geq j^*$.
Consequently, the $l$th probability in the last sum is bounded by
\begin{eqnarray} \label{c1}
&&\Pr\bigl( \|p_n(j^-) - p_n(l) -Ep_n(j^-) + Ep_n(l) \|_\infty\nonumber
\\[-8pt]\\[-8pt]
&&\qquad\quad\hspace*{-6pt} > T(n,j^-,l)
  + \bigl(\sqrt{0.99M}-2U(p_0,\Phi)\bigr) \sigma(l,n) \bigr)\nonumber
\end{eqnarray}
and we now apply Corollary \ref{est2} to this bound. Define the class
of functions
\[
\mathcal F := \mathcal F_{j^-,l} = \bigl\{2^{-l}\bigl(K_{j^-}(\cdot,y)-K_{l}(\cdot
,y)\bigr)/(4 \|\Phi\|_\infty) \bigr\},
\]
which is uniformly bounded by $1/2$ and satisfies (\ref{vcc}) for some
$A$ and $v$ independent of $l$ and $j^-$, by Lemma \ref{vc} (and a
computation on covering numbers). We compute $\sigma$, using (\ref
{secm}) and $l>j^-$:
\begin{eqnarray*}
\bigl(2^{-l}E(K_{j^-}-K_{l})(X,y)\bigr)^2 &\le& 2^{-2l+1}\bigl(EK_{j^-}^2(X,y) +
EK_l^2(X,y) \bigr)
\\
&\le& 2^{-2l+1} \|\Phi\|_2^2 \|p_0\|_\infty(2^{j^-} + 2^l) \le3 \cdot
2^{-l} \|\Phi\|_2^2 \|p_0\|_\infty,
\end{eqnarray*}
so that we can take $\sigma^2 = 3 \cdot2^{-l} \|\Phi\|_2^2 \|p_0\|
_\infty/(16 \|\Phi\|_\infty^2).$ The probability in (\ref{c1}) is then
equal to
\begin{eqnarray*}
&&\Pr\Biggl(\frac{2^{l} 4 \|\Phi\|_\infty}{n} \Biggl\|\sum_{i=1}^n f(X_i)-Pf \Biggr\|
_{\mathcal F}
\\
&&\qquad\quad\hspace*{-2pt}> \frac{2^l 4 \|\Phi\|_\infty}{n} 2\Biggl\|\sum_{i=1}^n
\varepsilon_i f(X_i)\Biggr\|_{\mathcal F}
+ \bigl(\sqrt{0.99M}-2U(p_0,\Phi
)\bigr)\sigma(l,n) \Biggr)
\\
&&\quad= \Pr\Biggl(\Biggl\|\sum_{i=1}^n f(X_i)-Pf \Biggr\|_{\mathcal F} > 2\Biggl\|\sum_{i=1}^n
\varepsilon_i f(X_i)\Biggr\|_{\mathcal F}
+ 3\frac{n (\sqrt
{0.99M}-2U(p_0,\Phi)) \sigma(l,n)}{3 \cdot2^l \cdot4\|\Phi\|_\infty} \Biggr).
\end{eqnarray*}
Since $n \sigma^2 / \log(1/\sigma) \simeq n/(2^ll) \to\infty$
uniformly in $l \in\mathcal J$, there exists $\lambda_n \to\infty$
independent of $l$ such that (\ref{sigmab}) is satisfied and the choice
\[
t =\frac{n (\sqrt{0.99M}-2U(p_0,\Phi)) \sigma(l,n)}{3 \cdot2^l \cdot
4\|\Phi\|_\infty}
\]
is admissible in Corollary \ref{est2} for $c_2(\lambda_n) = 1+120
\lambda_n^{-1} + 10{,} 800 \lambda_n^{-2}$. Hence, using Corollary \ref
{est2}, the last probability is bounded by
%
%e41 ###
\begin{equation} \label{finbd}
\le2 \exp\biggl(- \frac{n^2 (\sqrt{0.99M}-2U(p_0,\Phi))^2 (2^ll/n) 16 \|
\Phi\|_\infty^2 }{9 \cdot6.3 \cdot c_2(\lambda_n) 2^{2l} n 2^{-l} \|
\Phi\|_2^2 \|p_0\|_\infty16 \|\Phi\|_\infty^2} \biggr) \le2^{-l((q/2)+\delta)}
\end{equation}
for some $\delta>0$ and $q>1$, by the definition of $M$. Since $\sum_{l
\in\mathcal J: l \geq j} 2^{-l{(q/2)+\delta})} \le c 2^{-j((q/2)+\delta
)}$, we have proven (\ref{b0}).

\textit{Case II}: $\hat j_n = \tilde j_n$. The proof reduces to the
previous case since, by inequality (\ref{triang}), one has
\begin{eqnarray*}
&& \Pr(\{\tilde j_n^\varepsilon= j\} \cap\{0.99M \leq\tilde M\})
\nonumber
\\
&&\quad \leq \sum_{l \in\mathcal{J}: l \ge j} \Pr\bigl( \| p_n(j^-) - p_n(l) \|
_\infty> \bigl(B(\phi)+1\bigr) R(n,l) + \sqrt{0.99 M} \sigma(l,n)\bigr)
\\
&&\quad \le\sum_{l \in\mathcal{J}: l \ge j} \Pr\bigl( \| p_n(j^-) - p_n(l) \|
_\infty> T(n,j^-,l) + \sqrt{0.99 M} \sigma(l,n)\bigr).
\end{eqnarray*}
\upqed\end{pf*}

\section*{Acknowledgements} We would like to thank Patricia
Reynaud-Bouret and Benedikt P\"{o}tscher for helpful comments. The idea
of using Rademacher thresholds in Lepski's method arose from a
conversation with Patricia Reynaud-Bouret.

\printhistory

\end{document}